\nonstopmode

\documentclass[russian,english,11pt]{article}
\usepackage{babel}
\usepackage[T2A]{fontenc}
\usepackage[cp1251]{inputenc}
\usepackage{hyperref}
\usepackage{breakurl}
\usepackage{amssymb}
\usepackage{amsmath}
\usepackage{latexsym}
\usepackage{amscd}
\usepackage{amsfonts}
\usepackage{theorem}
\usepackage{mdwtab}
\usepackage{float}

\theoremstyle{change}

{\theorembodyfont{\rmfamily} \theoremstyle{change}

}

\newcommand{\ZZ}{\mathbb {Z}}
\newcommand{\RR}{\mathbb {R}}

\begin{document}

\title{Horizontal Monotonicity of the Modulus of the Riemann Zeta
Function  and  Related Functions}

\footnotetext[1] {2000 Mathematics Subject Classification : Primary 11M06, Secondary 11M26 \\
Keywords and phrases : Riemann zeta function, monotonicity}
\footnotetext[2] {The first named author was supported during this work by
grant НШ-3229.2012.1 of the President of the Russian Federation.}
\footnotetext[3] {The third named author was supported during
this work by a Discovery Grant from the Natural Sciences and
Engineering Research Council of Canada.}

\author{ Yuri Matiyasevich,  Filip Saidak,  Peter Zvengrowski }

\date{}

\maketitle

\begin{abstract}

\English
As usual let $s = \sigma + it$. For any fixed value $t = t_0$ with $|t_0| \geq 8$,
and for
$\sigma \leq 0$, we show that \ $|\zeta(s)|$ is strictly monotone decreasing
in $\sigma$, with the same result also holding for the related functions  $\xi$ of
Riemann and
 $\eta$ of Euler. The following inequality relating the
monotonicity of all three functions is proved:
$$ \Re\left(\frac {\eta'(s)}{\eta(s)} \right) <   \Re\left(\frac {\zeta'(s)}{\zeta(s)} \right)
         <   \Re\left(\frac {\xi'(s)}{\xi(s)} \right).
  $$
It is also shown that
extending the above monotonicity
result from $\sigma \leq 0$ to $\sigma \leq 1/2$, for any of $\zeta,\xi,\eta$, is equivalent to the
Riemann hypothesis.

\end{abstract}

\English

\section{Introduction} \label{sec:1}

Starting from the work of Riemann \cite{riemann}, the zeta function $\zeta(s)$
(with $s = \sigma + it$)    has
been primarily investigated in the vertical sense, especially in the
critical strip $0 \leq \sigma \leq 1$ and on the critical line $\sigma = 1/2$.
Questions related to the horizontal behaviour of $|\zeta(s)|$ (as usual we
write \ $s = \sigma + it$)  have been considered
by Saidak and Zvengrowski in \cite{saidak-zvengrowski}, and earlier by Spira
\cite {spira}. Indeed, the opening page of the article on the Riemann zeta
function in the Wolfram MathWorld \cite{wolfram} has a plot showing horizontal ``ridges" of
$|\zeta(\sigma + it)|$ for  $0 < \sigma < 1$ and $0 < t < 100$. To quote from this
article, ``the
fact that the ridges decrease monotonically for $0 < \sigma < 1/2$ is not a
coincidence since it turns out that monotone decrease implies the Riemann hypothesis,"
cf. \cite{saidak-zvengrowski} and \cite{borwein-bailey}. In this note, among other
things, we shall not only prove the converse of this assertion, but also the fact that
$|\zeta(s)|$ is monotone decreasing in $\sigma$ in the region $\sigma < 0$, subject
to the (minor) additional condition $|t| \geq 8$.

Recently a paper by Sondow and Dumitrescu \cite{sondow} has appeared
exploring this question for the related Riemann $\xi$ function (defined by
 $\xi(s): = (s-1)\Gamma(1 + s/2)
\pi^{-s/2}\zeta(s)$).
Here we shall
consider this question for $\zeta(s)$ (as mentioned above), as well as for
 $\xi(s)$ and Euler's function $\eta$ (cf. \cite{euler}) (also known as the
Dedekind $\eta$ function) defined by
$  \eta(s): = (1 - 2^{1-s})\zeta(s)$ or, for $\sigma > 0,$ \ by the alternating
Dirichlet series
$ \eta(s) = \sum_{n\geq 1} (-1)^{n+1}/n^s$.
A recent
paper of Srinivasan and Zvengrowski \cite{srinivasan-zvengrowski} also  examines
this question, for the $\Gamma$ function, and another recent paper of Alzer \cite{alzer},
titled ``Monotonicity Properties of the Riemann Zeta Function," concerns itself with
monotonicity of a function related to the zeta function, but only along the real line.
For completeness, let us
quote the results in \cite{sondow} and \cite{srinivasan-zvengrowski}.

{\bf   Theorem 1.1} (Sondow--Dumitrescu) : The xi function is increasing in modulus along every
horizontal half-line lying in any open right half-plane that contains no zeros of xi.
Similarly, the modulus decreases along each horizontal half-line lying in any zero-free,
open, left half-plane.

{\bf  Theorem 1.2} (Srinivasan--Zvengrowski) : For any fixed $t$ with $|t| > 5/4$, \ $|\Gamma(s)|$
\ is monotone increasing in $\sigma$.

 Section 2 starts by quoting an elementary lemma from \cite{srinivasan-zvengrowski} that
relates the horizontal increase or decrease of $|f(s)|$, for any holomorphic function
$f$, to \ $\Re(f'(s)/f(s))$, the real part of the logarithmic derivative of $f$. Using
this lemma we give a very short proof of the Sondow--Dumitrescu theorem. We also show
how a portion of this theorem was implicitly anticipated in a paper of P\'olya \cite
{polya} written in 1927. It is also related to work of Lagarias \cite{lagarias},
Haglund \cite{haglund}, and others, again this is briefly discussed in Section 2.

 In Section 3
we prove our first main result,  namely

{\bf   Theorem 1.3} : For \ $ |t| \geq 8, \ \ \sigma < \frac{1}{2}$, one has
$$ \Re\left(\frac {\eta'(s)}{\eta(s)} \right) <   \Re\left(\frac {\zeta'(s)}{\zeta(s)} \right)
         <   \Re\left(\frac {\xi'(s)}{\xi(s)} \right),
  $$
which relates the horizontal growth rates of all three functions under consideration. The
second main result, which follows as a corollary of this inequality together
with the results in Section 2, is now stated.

{\bf Theorem 1.4} : The moduli of all three functions $\eta(s), \zeta(s)$,
and $\xi(s)$ are monotone decreasing with respect to $\sigma$ in the region
\ $\sigma \leq 0$, \ $|t|\geq 8$. Extending this region to $\sigma \leq 1/2,$
for any of the three functions, is equivalent to the Riemann hypothesis.

The inequality given in Theorem 1.3 seems to indicate that in order to seek further results
on monotonicity, for $\sigma < 1/2$,  the most promising of the three functions
is $\eta$, and the least
promising $\xi$. On the other hand, combining the monotonicity results for $\zeta$
together with the Voronin Universality Theorem \cite{voronin} for $\zeta$
(or for $\log\zeta$) seems
to offer an approach to possibly showing that the Riemann hypothesis is false.  We also
note that the
inequality $|t|\geq 8$ is essential. Slightly
smaller numbers than $8$ will also work but, for $|t|<6.2897,$ the conclusion of
Theorem 3.4 is false for both $\zeta$
and $\eta$. Also,
for $\sigma > 1/2$, neither $|\zeta|$ nor $|\eta|$ are
monotone (by ``monotone" we always mean monotone with respect to $\sigma$).

\section{Monotonicity of $|\xi|$} \label{sec:2}

To measure the rate of change of $|f(s)|$ with respect to $\sigma$, the
following elementary lemma is useful. For
a proof cf. \cite{srinivasan-zvengrowski}.

{\bf  Lemma 2.1:} For any holomorphic function $f$, with $f(s)\not= 0$
in some open domain ${\cal D}$,
$$
\Re\left(\frac{f'(s)}{f(s)}\right) = \frac{1}{|f(s)|}
\cdot \frac{\partial |f(s)|}{\partial \sigma}, \ \ \  s \in {\cal D} \ .  $$

{\bf  Corollary 2.2:} For \ $s\in {\cal D}, \ \ \ {\displaystyle {\rm sgn}
\left(\frac{\partial |f(s)|}
{\partial \sigma}\right) = {\rm sgn}\left( \Re\left(\frac{f'(s)}{f(s)}\right) \right)}$.

The fact that Lemma 2.1 does not apply at a zero of $f$ is not a problem towards
our main objectives, as the next lemma shows.

{\bf  Lemma 2.3:} (a) \  Let $f$ be holomorphic in an open domain ${\cal D}$ and not
identically zero. Let us also suppose
$\Re(f'(s)/f(s))< 0$  for all \ $s \in {\cal D}$ such that
\ $f(s)\not= 0.$ Then $|f(s)|$ is strictly decreasing  with respect
to $\sigma$ in ${\cal D}$, i.e. for each \ $s_0 \in {\cal D}$ \ there exists a
$\delta > 0$ such that $|f(s)|$ is strictly monotonically decreasing with respect to
$\sigma$ on the horizontal interval from \ $s_0 - \delta$\ to \ $s_0 + \delta$.

(b) \ \ Conversely, if $|f(s)|$ is decreasing
with respect to $\sigma$ in ${\cal D}$, then
$\Re(f'(s)/f(s)) \leq 0$  for all \ $s \in {\cal D}$ \ such that
 \ $f(s) \not= 0$.

Proof of (a): From Lemma 2.1 and Corollary 2.2 it clearly suffices to show this
for those $s_0 = \sigma_0 + it_0 \in {\cal D},$ \ where \ $f(s_0) = 0$. Thanks to
$f$ being holomorphic and not identically $0$ there exists $\delta > 0$
with \ $\{ s : |s-s_0| < \delta\}
\subset {\cal D}$ and with no further zeros of $f$ in this open disc.   Then
using the next part of the hypothesis and Corollary 2.2, \ $|f(s)|$ is strictly
decreasing with respect to $\sigma$ on the two horizontal open intervals from \
$\sigma_0 - \delta + it_0$ \ to $\sigma_0 + it_0$, and from \  $\sigma_0 + it_0$
\ to  $\sigma_0 + \delta + it_0$. Since $|f|$ is continuous in ${\cal D}$, a
simple continuity argument shows that it must be strictly decreasing on the
entire horizontal interval from \ $\sigma_0 - \delta + it_0$ \ to \
 $\sigma_0 + \delta + it_0$.

Proof of (b):  Conversely, we are assuming \ $\partial|f(s)|/\partial\sigma \leq 0$ \
in ${\cal D}$, so
Lemma 2.1 implies that $\Re(f'(s)/f(s))\leq  0$ at any $s\in {\cal D}$ for which
$f(s) \not= 0$.   \hfill $\Box$

Of course the analogous results hold for monotone increasing and \linebreak
$\Re(f'(s)/f(s))> 0$. Combining Lemma 2.3 with the fact that a function can have
no zeros in an open domain in which its modulus is stictly monotone decreasing
(increasing) with respect to $\sigma$ gives the next result.

{\bf  Corollary 2.4:} With the same hypotheses as in Lemma 2.3 (a), $f$ has no
zeros in ${\cal D}$.

Let us now apply the above to the Riemann $\xi$ function and thereby
give a short proof of Theorem 1.1. It is well known
that $\xi(1-s) =
\xi(s)$ and that \ $\xi(\overline s) = \overline{ \xi(s)}$  . Hence
$|\xi(1/2 - \sigma + it) | = |\xi(1/2 + \sigma -it)| =
|\xi(1/2 + \sigma +it)|   $, which shows that $|\xi|$  is
symmetrical about the critical line $\sigma = 1/2.$ So showing that $|\xi|$ is
monotone decreasing in a domain to the left of the critical line is
equivalent to
showing it is monotone increasing  in the
reflexion of the  same domain about the point $s = 1/2$, and
this is what we shall show.

{\bf  Theorem 1.1:}
(Sondow--Dumitrescu) Let $\sigma_0$ be greater than  or
equal to the
real part of any zero of $\xi$. Then $|\xi(s)|$ is strictly monotone
increasing in the half plane $\sigma > \sigma_0$.

Proof: We start with the formula due to  Hadamard \cite{hadamard} and
von Mangoldt \cite{vonMangoldt} (also  cf.
\cite{ksz}, (36), or simply take the logarithmic derivative of the final formula
given in \cite{edwards}, \S 2.8)
$$ \frac{\xi'(s)}{\xi(s)} = \sum_{\rho} \frac{1}{s - \rho},$$
where the summation is taken over all zeros $\rho$ of $\xi$
(which, as is well known, lie in the critical strip $0 < \Re(\rho) < 1$),
in conjugate
pairs and in order of increasing $\Im(\rho)$.
If any such zero be written as \ $\rho = \alpha + i\beta,$ then by
hypothesis \ $\sigma > \alpha$. It is then trivial to check that
\ $\Re(1/(s - \rho)) = (\sigma - \alpha)/[(\sigma - \alpha)^2 +
(t - \beta)^2]  > 0$,\ hence \ $\Re(\xi'(s)/\xi(s)) > 0$ \ and by
Corollary 2.4 $|\xi(s)|$ is monotone increasing in $\sigma$, in the
given half plane $\sigma > \sigma_0.$ \hfill $\Box$

Combining this theorem with well known facts about the zeros of $\xi$, and the fact
that a function can have no zeros in an open domain where its modulus is monotone
increasing (decreasing), gives the next result.

{\bf  Corollary 2.5:} (Sondow--Dumitrescu) In the right (left) half plane
  $\sigma \geq   1  \ (\sigma \leq 0$),
 \ $|\xi|$ is monotone increasing (decreasing). The same is true for the right (left)
half plane
 $\sigma \geq 1/2$ ($\sigma \leq 1/2$) if and only if the Riemann hypothesis is true.

The second part of this corollary, which is the same as Corollary 1 in
\cite{sondow}, is actually implicit in a paper written by P\'olya in
1927 \cite{polya} which discusses the ``Nachlass" of J.L.W.V. Jensen, after
suitable interpretation. Namely, following ${\rm I}'$ on p. 18 of \cite{polya}, and
using $z = x + iy$ as in this reference, we consider the holomorphic
function \ $F(z) = \xi(1/2 - iz) = \xi(1/2 + y - ix)$. Note that
$|F(z)| = |\xi(1/2 + y - ix)| = |\xi(1/2 + y + ix)|$, since
 \ $\xi(\overline{s}) = \overline{\xi(s)}.$ The condition that
all zeros of $F$ are real is precisely the Riemann hypothesis, indeed this
was Riemann's original formulation. According to condition ${\rm I}'$, this is
equivalent to \ ${\displaystyle \frac{\partial^2|\xi(1/2 + y + ix)|^2}
{\partial^2 y} \geq 0}.$ \ This implies that  \ $|\xi(1/2 + y + ix)|^2$ \
is a convex function of $y$. By symmetry it has  zero derivative at
$y = 0$, hence it is monotone increasing for $y \geq 0$ and monotone
decreasing for $y \leq 0$. The same is then also true for
$|\xi(1/2 + y + ix)|$.   And conversely, as already remarked before
Corollary 2.5, these monotonicity properties imply the Riemann hypothesis.

The fact that \ $ \Re(\xi'(s)/\xi(s)) > 0$ \ when $\sigma > 1$, and that the
Riemann hypothesis is equivalent to the same statement for \ $\sigma > 1/2$
(for which we gave a short proof above) also
appears in the 1999 paper of Lagarias \cite{lagarias} and the 1997 paper of
Hinkkanen \cite{hink}. Combining this with Lemma 2.3 gives an immediate proof of
Theorem 1.1. Another version  the Sondow-Dumitrescu result appears as a ``known result"
at the beginning of Section 6 of \cite{haglund}, this time for the related $\Xi$
function
(the horizontal monotonicity of $\xi$ being equivalent to vertical monotonicity of
$\Xi$), but no reference or proof is given.

\section{Proof of Theorems 1.3 and 1.4  } \label{sec:3}

For convenience we label the first inequality of Theorem 1.3 as (A), and the second (B).
To prove either of these we shall take the logarithmic
derivatives of the formulae given for $\xi, \eta$ in the Introduction, and then look
at the real part of these
logarithmic derivatives.  Again, for convenience, we will divide the proof into
corresponding  parts (A), (B), and separately give two lemmas that will be of use.

{\bf  Lemma 3.1:} For $\sigma < 1$, one has
${\displaystyle \Re\left(\frac{1}{2^{s-1}-1}\right) < 0}$.

Proof: First note that $2^{s-1} - 1 = 0$ if and only if \ $s = 1 + 2n\pi i/\log 2, \ \
n\in \ZZ.$ In particular $2^{s-1} - 1 \not= 0$ for $\sigma < 1$. Now
$$\Re\left(\frac{1}{2^{s-1}-1}\right) = \frac{2^{\sigma - 1} \cos(t\log 2) - 1}
{|2^{s-1} - 1 |^2}.$$
The denominator of this expression is strictly positive since $\sigma < 1$. As for the
numerator, one has \ \ $|2^{\sigma - 1} \cos(t\log 2)| < |\cos(t\log 2)| \leq 1,     $
so the numerator is strictly negative. \hfill $\Box$

Proof of (A): \ From the formula given at the beginning of the
Introduction for $\eta(s)$, it follows that \  $\log(\eta(s))
 = \log(1 - 2^{s-1}) + \log(\zeta(s))$. Differentiating,
$$\frac{\eta'(s)}{\eta(s)} = \frac{2^{1-s} \log 2}{1 - 2^{1-s}} + \frac {\zeta'(s)}{\zeta(s)}
 = \frac{\log 2}{2^{s-1} - 1} + \frac {\zeta'(s)}{\zeta(s)}.$$
Taking the real parts, and using Lemma 3.1 as well as $\log 2 > 0$, completes the proof
(indeed for \ $\sigma < 1$).

For the second inequality it will be necessary to recall the digamma function \ $\Psi(s):
= \Gamma'(s)/\Gamma(s)$. We list a few of its properties as the next lemma.

{\bf  Lemma 3.2:}\ (i) $\Psi(s) - \Psi(1-s) = -\pi \cot(\pi s)$, \\
\smallskip
 (ii) \ \ $|\Re(\Psi(s)) - \Re(\Psi(1-s))| < 3\pi e^{-2\pi t}, \ \ t \geq 0.1   $,\\
\smallskip
(iii) \ \ in the sector of the complex plane \ $-\theta < \arg(s) < \theta$, one has
$$\Psi(s) = \log(s) - \frac{1}{2s} + R_0'(s), \ \ {\rm where} \ \
|R_0'(s)| \leq \sec^3(\theta/2)\cdot \frac{B_2}{2|s|^2},$$
with $B_2 = 1/6$ being the second Bernoulli number,\\
(iv) \ \ $|x/(x^2 + t^2)| \leq 1/(2|t|)$, for any \ $x,t \in \RR, \ t \not= 0$,\\
(v) \ for any $\sigma \geq 1/2$, and $|t| \geq 8$,  \ $\Re(\Psi(s)) > 2.0096.$

Formula (i) is a simple consequence of Euler's reflexion formula for the $\Gamma$ function,
it can be found e.g. in \cite{srivastava-choi}, p.14. Formula (ii) follows from (i) and
doing an elementary estimate of $\Re(\cot(z))$, since (i) implies
$$|\Re(\Psi(s))- \Re(\Psi(1-s))| \leq |\Psi(s)- \Psi(1-s)| = \pi|\cot(z)|,$$
where for convenience we set \ $\pi s =  z = x + iy $.
We outline the remaining  steps towards proving (ii), which are
essentially an exercise in calculus.
First recall that
$$ \cot(z) = \frac{\cos x  \cosh y - i \sin x \sinh y}{\sin x \cosh y
+ i \cos x \sinh y}.$$
From this it is easy to derive
$$\Re(\cot(z)) = \frac{\sin(2x)}{b - \cos(2x) + 1} =: g_b(x),
 \ {\rm where} \ b = 2\sinh^2(y) > 0.$$
We claim that \ $|g_b(x)| < 3e^{-2y},$   \ when $y > \log(3)/4$. Indeed, using elementary
calculus one shows that \ ${\displaystyle |g_b(x)|_{\rm max} = \frac{1}{\sqrt{b^2 + 2b}}}$,
hence proving the claim reduces to showing \
${\displaystyle  \frac{1}{\sqrt{b^2 + 2b}}< 3e^{-2y}   }$. Using the definition of
$b$, this inequality reduces to \ $y > \log(3)/4$ and the claim is thus  proved.
Finally, substituting $z = \pi s$, we obtain (ii) with \ $y = \pi t > \log(3)/4$, \ i.e.
\ $t > \log(3)/(4\pi) = .08742..$ \ .

Formula (iii) is a special case ($n = 0$) of
the Stirling series
$$\Psi(s) = \log(s) - \frac{1}{2s} -\sum_{k=1}^n
\frac{B_{2k}}{2k s^{2k}}  + R'_{2n} \ $$
for digamma, together with the Stieltjes estimate for
the error term (cf. \cite{edwards}, p.114,  or
the original manuscript
of Stieltjes \cite{stieltjes})
$$|R'_{2n}| \leq (\sec(\theta/2))^{2n+3} \left|\frac{B_{2n+2}}{(2n+2)s^{2n+2}}
\right|.$$
Formula (iv) is equivalent to \ $0 \leq (|x| - |t|)^2$. Finally, from
(iii) applied to the sector $-\pi/2 < \theta < \pi/2$, we have
$$\Re(\Psi(s)) = \log|s| -\frac{\sigma}{2|s|^2} + \Re(R_0'(s)),$$
where \ \ $|\Re(R_0'(s))| \leq |R_0'(s)| < 2\sqrt{2}/(6|s|^2)$. Now assume as in (v)
that $\sigma \geq 1/2$ and $|t|\geq 8$, then using this estimate for the remainder
term as well as $|s| > 8$, we obtain
$$\Re(\Psi(s))  \geq \log 8 - 1/16 - \sqrt{2}/(3\cdot 64) = 2.0096, $$
where (iv) was used to give the
$1/16$ estimate for the second term. This completes the proof of Lemma 3.2.

Proof of (B): \ Taking the logarithmic derivative of the formula given at the beginning
of the Introduction for $\xi(s)$ gives
$$\frac{\xi'(s)}{\xi(s)} =  \frac{\zeta'(s)}{\zeta(s)} + \frac{1}{s-1}
    +\frac{1}{2}\Psi\left(\frac {s}{2} + 1\right) - \frac{1}{2}\log \pi.$$
Hence, to complete the proof of (B), it suffices to show that
$$\Re\left(\frac{1}{s-1} + \frac{1}{2}\Psi\left(\frac {s}{2} + 1\right)\right)-
\frac{1}{2}\log \pi > 0, \ \ \sigma < \frac{1}{2}, \ \ 8 \leq |t|.$$
Now, by Lemma 3.2 (iv), the first term is greater than or equal to $-1/16$. By
Lemma 3.2 (v) \ $(1/2)\Re(\Psi(z)) > 1.0048$, \ at least when $\sigma \geq 1/2$ and
 \ $|t| \geq 8$.
However, applying Lemma 3.2 (ii), we see that the same holds for any $z$ with $|t|
\geq 8$, at least to within \ $3\pi e^{-16\pi}$ \ which is negligible here. Thus the
sum in question is greater than \ $- 1/16 + 1.0048 - \log \pi/2 > 0.$ \hfill $\Box$

{\bf Theorem 3.4:} The moduli of all three functions $\eta(s), \zeta(s)$,
and $\xi(s)$ are monotone decreasing with respect to $\sigma$ in the region
\ $\sigma \leq 0$, \ $|t|\geq 8$. Extending this region to $\sigma \leq 1/2,$
for any of the three functions, is equivalent to the Riemann hypothesis.

Proof: For $\sigma \leq 0$, we have seen in the proof of Theorem
2.5 that $\Re(\xi'(s)/\xi(s)) < 0$. Combining this with the inequalities in
Theorem 3.1 shows that the same is true for $\zeta$ and $\eta$, thus all
three are monotone decreasing in modulus for \ $\sigma \leq 0, \ |t| \geq 8$. And
the same argument used in Corollary 2.5 shows that extending this to the larger
region \ $\sigma \leq 1/2, \ |t| \geq 8$, is equivalent to the Riemann hypothesis.
\hfill $\Box$

\English

\noindent
\author{Yuri Matiyasevich}\\
{\small St.Petersburg Department}\\
{\small \phantom{xx}of Steklov Institute of Mathematics}\\
{\small  \phantom{xx}of Russian Academy of Sciences}\\
{\small  \phantom{xx}(POMI RAN)}\\
{\small 27, Fontanka } \\
\small {St. Petersburg, 191023, Russia} \\
%\small {} \\
\small {e-mail: \url{yumat@pdmi.ras.ru}}\\

\noindent
\author{Filip Saidak}\\
{\small Department of Mathematics } \\
\small {University of North Carolina} \\
\small {Greensboro, NC 27402, U.S.A.} \\
\small {e-mail: \url{saidak@nz11.com }}\\

\noindent
\author{Peter Zvengrowski} \\
{\small {Department of Mathematics and Statistics} \\
\small{University of Calgary} \\
\small {Calgary, Alberta, Canada T2N 1N4}\\
\small { e-mail: \url{zvengrow@ucalgary.ca}}


\begin{thebibliography}{99}


\bibitem{alzer} {Alzer, H.}{\it Monotonicity properties of the Riemann zeta function},
Meditterean J. Math (2011)

\bibitem{borwein-bailey}{Borwein, J., Bailey, D.}, \ Mathematics by Experiment -- Plausible
Reasoning in the Twenty-first Century, A.K.Peters, Wellesley, Massachusetts (2003).

\bibitem{edwards}{Edwards, H.M.}, \ Riemann's Zeta Function, Academic Press, New York (1974).
Reprinted by Dover Publications, Mineola, N.Y. (2001).

\bibitem{euler}{Euler, L.}, \ {\it Remarques sur un beau rapport entre les s\'eries des
puissances tant directes que r\'eciproques} (lu en 1749) (presented 1761), Hist.
Acad. Roy. Sci. Belles-Lettres Berlin {\bf 17} (1768), 83--106. Also in ``Opera Omni,"
Ser. 1, Vol. 15, 70--90. Original and English translation available via
\url{http://www.math.dartmouth.edu/~euler/pages/E352.html}


\bibitem{hadamard}{Hadamard, J.}, \ {\it \'Etude sur les propri\'et\'es des fonctions ent\^eires
et en particulier d'une fonction consid\'er\'ee par Riemann}, J. Math. Pures Appl.
(4) {\bf 9} (1893), 171--215.


\bibitem{haglund}{Haglund, J.}, \ {\it Some conjectures on the zeros of approximates to the
Riemann $\Xi$ function and incomplete gamma functions,} Central European J. Math. {\bf 9}
(2) (2011), 302--318.

\bibitem{hink}{Hinkkanen, A.}\ {\it On functions of bounded type}, Complex Variables
{\bf 34} (1997), 119--139.

\bibitem{ksz}{Kudryavtseva, E., Saidak, F., Zvengrowski, P.},\ {\it Riemann and his zeta function}, \ Morfismos {\bf 9},
No. 2 (2005), 1--48.

\bibitem{lagarias}{Lagarias, J.C.},{\it On a positivity property of the Riemann $\xi$ function}, Acta Arith. {\bf 89} (3) (1999), 217--234.

\bibitem{vonMangoldt}{von Mangoldt, H.}, \ {\it Zu Riemann's Abhandlung `Ueber die Anzahl der Primzahlen
unter einer gegebenen Gr\"osse'}, J. Reine Angew. Math. {\bf 114} (1895), 255--305.


\bibitem{murty}{Murty, M. R.}, \  Problems in Analytic Number Theory, \ Graduate Texts in Math. {\bf 206}, Springer-Verlag,
N.Y., Berlin, Heidelberg (2001).

\bibitem{polya}{P\'olya, G.}, \ {\it \"Uber die algebraisch-funktiontheoretischen
Untersuchungen von J.L.W.V. Jensen}, Kgl. Danske Vid. Sel. Math.-Fys.
Medd. {\bf 7} (1927), 3--33. Also available in P\'olya, George, Collected Papers Vol. II,
Location of Zeros, edited by R. P. Boas, Mathematicians of Our Time, Vol. 8, The MIT
Press, Cambridge, Mass., London (1974).

\bibitem{riemann}{Riemann, B.}, \ {\it \"Uber die Anzahl der Primzahlen
unter einer gegebenen Gr\"osse} (1859, \emph{Monatsberichter der Berliner Akademie},
November 1859. Included in : Riemann, B., Gesammelte Werke,
 Teubner,
Leipzig, 1892; reprinted by Dover Books, New York, 1953. Original and English
translation available via \url{http://www.claymath.org/millenium/Riemann_Hypothesis/
1859_manuscript/}

\bibitem{saidak-zvengrowski}{Saidak, F., Zvengrowski, P.}, \ {\it On the modulus of the Riemann zeta function
in the critical strip}, \ Math. Slovaca {\bf 53} No. 2 (2003), 145--172.



\bibitem{sondow}{Sondow, J., Dumitrescu, C.}, \ {\it A monotonicity property
of Riemann's xi function and a reformulation of the Riemann hypothesis},
\ Periodica Math. Hung. {\bf 60} I (2010), 37--40. Also available at 
\url{http://arxiv.org/abs/1005.1104}. 

\bibitem{spira}{Spira, R.}, \ {\it An inequality for the Riemann zeta-function}, \ Duke Math. J. {\bf 32}
(1965), 247--250.


\bibitem{srinivasan-zvengrowski}{Srinavasan, G.K. , Zvengrowski, P.}, \ {\it On the horizontal
monotonicity of $|\Gamma(s)|$}, \ Can. Math. Bull.{\bf 54} (3) (2011), 538--543.

\bibitem{srivastava-choi}{Srivastava, H.M., Choi, J.}, Series Associated
with the Zeta and Related Functions, \ Kluwer
Academic Publishers, Dordrecht (2001).

\bibitem{stieltjes} {Stieltjes, T. J.}, \ {\it Sur le developpement de {\rm log}$\Gamma(a)$,}
J. Math. Pures Appl.(9) {\bf 5} (1889), 425--444.

\bibitem{voronin} \Russian Воронин, С.\,М. [{Voronin, S. M.}],
\emph{ Теорема об ``универсальности'' дзета-функции Римана}
\ [{\it Theorem on the universality of the Riemann zeta
function}], Изв. Акад. Наук СССР, сер. матем.  {\bf 39:3} (1975),  475–486.
\English Translated in
Math. USSR Izv. {\bf 9} (1975), 443--445.

\bibitem{wolfram} \ Wolfram MathWorld \ \url{http://
mathworld.wolfram.com/RiemannZetaFunction.html}



\end{thebibliography}
\end{document}